\documentclass[12pt]{article}
\input{amssym}
\newtheorem{prethm}{{\bf Theorem}}

\newenvironment{thm}{\begin{prethm}{\hspace{-0.5
em}{\bf.}}}{\end{prethm}}
\newtheorem{precor}{{\bf Corollary}}

\newenvironment{cor}{\begin{precor}{\hspace{-0.5
em}{\bf.}}}{\end{precor}}
\newtheorem{preprop}{{\bf Proposition}}

\newtheorem{preque}{{\bf Question}}

\newtheorem{preques}{{\bf Question}}

\newtheorem{prelemma}{{\bf Lemma}}

\newenvironment{lemma}{\begin{prelemma}{\hspace{-0.5
em}{\bf.}}}{\end{prelemma}}
\newtheorem{prefact}{{\bf Fact}}

\newtheorem{preobs}{{\bf Observation}}

\newenvironment{obs}{\begin{preobs}{\hspace{-0.5
em}{\bf.}}}{\end{preobs}}
\newtheorem{prefig}{{\bf Figure}}

\newtheorem{prelemm}{{\bf Lemma}}

\newtheorem{preex}{{\bf Example}}

\newtheorem{prepro}{{\bf Proposition}}

\newtheorem{prelem}{{\bf Theorem}}

\newenvironment{lem}{\begin{prelem}{\hspace{-0.5
em}{\bf.}}}{\end{prelem}}
\newtheorem{preproof}{{\bf Proof.}}

\newenvironment{proof}[1]{\begin{preproof}{\rm
               #1}\hfill{$\rule{2mm}{2mm}$}}{\end{preproof}}

\newtheorem{preconj}{{\bf Conjecture}}

\newtheorem{predeff}{{\bf Definition}}

\def\newpic#1{}
\date{}
\begin{document}

\title{
{\Large{\bf On Randomly $\bf k$-Dimensional Graphs}}}
%

{\small
\author{
{\sc Mohsen Jannesari} and {\sc Behnaz Omoomi  }\\
[1mm]
{\small \it  Department of Mathematical Sciences}\\
{\small \it  Isfahan University of Technology} \\
{\small \it 84156-83111, Isfahan, Iran}}




 \maketitle \baselineskip15truept

\begin{abstract}
For an ordered set $W=\{w_1,w_2,\ldots,w_k\}$ of vertices and a
vertex $v$ in a connected graph $G$, the ordered $k$-vector
$r(v|W):=(d(v,w_1),d(v,w_2),\ldots,d(v,w_k))$ is  called  the
(metric) representation of $v$ with respect to $W$, where $d(x,y)$
is the distance between the vertices $x$ and $y$. The set $W$ is
called  a resolving set for $G$ if distinct vertices of $G$ have
distinct representations with respect to $W$. A resolving set for
$G$ with minimum cardinality is called a basis of $G$ and its
cardinality is the metric dimension of $G$.
 A connected graph $G$ is called randomly $k$-dimensional
graph if each $k$-set of vertices of $G$ is a basis of $G$. In
this paper, we study randomly $k$-dimensional graphs and provide
some properties of these graphs.
\end{abstract}

\section{Introduction}
We refer to~\cite{west} for graphical notations and terminologies
not described in this paper. Throughout the paper, $G=(V,E)$ is a
finite, simple, and connected  graph. The distance between two
vertices $u$ and $v$, denoted by $d(u,v)$, is the length of a
shortest path between $u$ and $v$ in $G$. Also, $N(v)$ is the set
of all neighbors of vertex $v$ and $deg(v)=|N(v)|$ is the degree
of vertex $v$. The maximum degree of the graph $G$, $\Delta(G)$ is
$max_{v\in V(G)}deg(v)$. We mean by $\omega(G)$, the number of
vertices in a maximum clique in $G$. For a subset $S$ of $V(G)$,
$G\setminus S$ is the induced subgraph $\langle V(G)\setminus
S\rangle$ by $V(G)\setminus S$ of $G$. A set $S\subseteq V(G)$ is
a separating set in $G$ if $G\setminus S$ has at least two
connected components. We call a vertex $v\in V(G)$ a cut vertex of
$G$ if $\{v\}$ is a separating set in $G$. If $G\neq K_n$ has no
cut vertex, then $G$ is called a $2$-connected graph. The
notations $u\sim v$ and $u\nsim v$ denote the adjacency and
non-adjacency relation between $u$ and $v$, respectively. The
symbol $(v_1,v_2,\ldots, v_n)$ represents a path of order $n$,
$P_n$.

For an ordered set $W=\{w_1,w_2,\ldots,w_k\}\subseteq V(G)$ and a
vertex $v$ of $G$, the  $k$-vector
$$r(v|W):=(d(v,w_1),d(v,w_2),\ldots,d(v,w_k))$$
is called  the ({\it metric}) {\it representation}  of $v$ with
respect to $W$. The set $W$ is called a {\it resolving set} for
$G$ if distinct vertices have different representations.  A
resolving set for $G$ with minimum cardinality is  called a {\it
basis} of $G$, and its cardinality is the {\it metric dimension}
of $G$, denoted by $\beta(G)$.

 For example, the graphs $G$ and $H$
in Figure~\ref{fig} have the basis $B=\{v_1,v_2\}$ and hence
$\beta(G)=\beta(H)=2$. The representations of vertices of $G$ with
respect to $B$ are
$$r(v_1|B)=(0,1),\quad r(v_2|B)=(1,0),\quad r(v_3|B)=(2,1),\quad
r(v_4|B)=(2,2),\quad r(v_5|B)=(1,2).$$ Also, the representations
of vertices of $H$ with respect to $B$ are
$$r(v_1|B)=(0,1),\quad r(v_2|B)=(1,0),\quad r(v_3|B)=(1,1),\quad
r(v_4|B)=(2,2),\quad r(v_5|B)=(1,2).$$
\begin{figure}[ht]
\centering \unitlength=.4mm
\caption{\label{fig} \scriptsize{$bas(G)=\beta(G)=res(G)$ and
$bas(H)\neq\beta(H)\neq res(H)$.}}
\end{figure}
\\To see that whether a given set $W$ is a resolving set for $G$,
it is sufficient to look at the representations of vertices in
$V(G)\backslash W$, because $w\in W$ is the unique vertex of $G$
for which $d(w,w)=0$. When $W$ is a resolving set for $G$, we say
that $W$ {\it resolves} $G$. In general, we say an ordered set $W$
resolves a set $T$ of vertices in $G$, if the representations of
vertices in $T$ are distinct with respect to $W$. When $W=\{x\}$,
we say that  vertex $x$ resolves $T$.
\par In~\cite{Slater1975}, Slater introduced the idea of a resolving
set and used a {\it locating set} and the {\it location number}
for what we call a resolving set and the metric dimension,
respectively. He described the usefulness of these concepts when
working with U.S. Sonar and Coast Guard Loran stations.
Independently, Harary and Melter~\cite{Harary} discovered the
concept of the location number as well and called it the metric
dimension. For more results related to these concepts
see~\cite{cartesian product,bounds,sur1,landmarks,sur2}. The
concept of a resolving set has various applications in diverse
areas including coin weighing problems~\cite{coin}, network
discovery and verification~\cite{net2}, robot
navigation~\cite{landmarks}, mastermind game~\cite{cartesian
product}, problems of pattern recognition and image
processing~\cite{digital}, and combinatorial search and
optimization~\cite{coin}.
\par The following simple result is very useful.
\begin{obs}~\rm\cite{extermal}\label{twins}
Let $G$ be a graph and $u,v\in V(G)$ such that,
$N(v)\backslash\{u\}=N(u)\backslash\{v\}$. If $W$ resolves $G$,
then $u$ or $v$ is in $W$.
\end{obs}
It is obvious that for a graph $G$ of order $n$,
$1\leq\beta(G)\leq n-1$.
\begin{lem}~\rm\cite{Ollerman}\label{B=1,B=n-1}
Let $G$ be a graph of order $n$. Then,\begin{description}\item (i)
$\beta(G)=1$ if and only if $G=P_n$,\item (ii) $\beta(G)=n-1$ if
and only if $G=K_n$.
\end{description}
\end{lem}
\par The {\it basis number}, $bas(G)$, of $G$ is the maximum integer $r$
such that, every $r$-set of vertices of $G$ is a subset of some
basis of $G$. Also, the {\it resolving number}, $res(G)$, of $G$
is the minimum integer $k$ such that, every $k$-set of vertices of
$G$ is a resolving set for $G$. These parameters are introduced
in~\cite{basis} and ~\cite{res(G)}, respectively. Clearly, if $G$
is a graph of order $n$, then  $0\leq bas(G)\leq \beta(G)$ and
$\beta(G)\leq res(G)\leq n-1$. Chartrand et al. in~\cite{basis}
considered graphs $G$ with $bas(G)=\beta(G)$. They called these
graphs {\it randomly $k$-dimensional} graphs, where $k=\beta(G)$.
Obviously, $bas(G)=\beta(G)$ if and only if $res(G)=\beta(G)$. In
the other word, a randomly $k$-dimensional graph is a graph which
every $k$-set of its vertices is a basis. For example in graph $G$
of Figure~\ref{fig}, if $W$ is a set of two adjacent vertices,
then the representations of vertices in $V(G)\setminus W$ with
respect to $W$ are $(1,2),~(2,2)$, and $(2,1)$. Also, if $W$ is a
set of two non-adjacent vertices, then the representations of
vertices in $V(G)\setminus W$ with respect to $W$ are
$(1,1),~(1,2)$, and $(2,1)$. Therefore, $G$ is a randomly
$2$-dimensional graph. But, in graph $H$ of Figure~\ref{fig},
$\{v_1,v_4\}$ is not a resolving set, hence $H$ is not a randomly
$2$-dimensional graph. Since $\{v_1,v_2\}$, $\{v_1,v_3\}$, and
$\{v_4,v_5\}$ are bases of $H$, $bas(H)=1$. Also, $res(H)=3$,
because every $3$-set of $V(H)$ is a resolving set in $H$.
\par Obviously, $K_1$ and $K_2$ are the only
randomly $1$-dimensional graphs. Chartrand et al.~\cite{basis}
proved that a graph $G$ is randomly $2$-dimensional if and only if
$G$ is an odd cycle. In this paper, we first characterize all
graphs of order $n$ and resolving number $1$ and
 $n-1$. Then, we provide
some properties of randomly $k$-dimensional graphs.
\section{Main Results}\label{main}
We first characterize all graphs $G$ with $res(G)=1$ and all
graphs $G$ of order $n$ with $res(G)=n-1$.
\begin{thm} \label{res=1,res=n-1}
Let $G$ be a graph of order $n$. Then,
\begin{description} \item (i) $res(G)=1$ if and only if
$G\in\{P_1,P_2\}$.\item (ii) $res(G)=n-1$ if and only if
$N(v)\backslash\{u\}=N(u)\backslash\{v\}$, for some $u,v\in V(G)$.
\end{description}
\end{thm}
\begin{proof}{ (i) It is easy to see that for $G\in
\{P_1,P_2\}$, $res(G)=1$. Conversely, let $res(G)=1$. Thus, $1\leq
\beta(G)\leq res(G)=1$ and hence, $\beta(G)=1$. Therefore, by
Theorem~\ref{B=1,B=n-1}, $G=P_n$. If $n\geq3$, then $P_n$ has a
vertex of degree $2$ and this vertex does not resolve its
neighbors. Thus, $res(G)\geq 2$, which is a contradiction.
Consequently, $n\leq2$, that is $G\in\{P_1,P_2\}$.\\
(ii) Let $u,v\in V(G)$ such that,
$N(v)\backslash\{u\}=N(u)\backslash\{v\}$. If $res(G)\leq n-2$,
then the set $V(G)\setminus \{u,v\}$ is a resolving set for $G$.
But, by Observation~\ref{twins}, every resolving set for $G$
contains at least one of the vertices $u$ and $v$. This
contradiction implies that, $res(G)=n-1$. Conversely, let
$res(G)=n-1$. Thus, there exists a subset $T$ of $V(G)$ with
cardinality $n-2$ such that, $T$ is not a resolving set for $G$.
Assume that, $T=V(G)\setminus \{u,v\}$. If
$N(u)\backslash\{v\}\neq N(v)\backslash\{u\}$, then there exists a
vertex $w\in T$ which is adjacent to only one of the vertices $u$
and $v$ and hence, $d(u,w)\neq d(v,w)$. Since $w\in T$, $T$
resolves $G$, which is a contradiction. Therefore,
$N(u)\backslash\{v\}=N(v)\backslash\{u\}$.}\end{proof}
\begin{cor}\label{not twin} If $G\neq K_n$ is a randomly
$k$-dimensional graph, then  for each pair of vertices $u,v\in
V(G)$, $N(v)\backslash\{u\}\neq N(u)\backslash\{v\}$.
\end{cor}
\begin{proof}{
If $N(v)\backslash\{u\}=N(u)\backslash\{v\}$ for some $u,v\in
V(G)$, then by Theorem~\ref{res=1,res=n-1}, $res(G)=n-1$, where
$n$ is the order of $G$. Since $G$ is a randomly $k$-dimensional
graph, $\beta(G)=res(G)=n-1$. Therefore, by
Theorem~\ref{B=1,B=n-1}, $G=K_n$, which is a contradiction. Hence,
for each $u,v\in V(G)$, $N(v)\backslash\{u\}\neq
N(u)\backslash\{v\}$. }\end{proof}
\begin{lemma}\label{delta>2} If $G$ is a randomly
$k$-dimensional graph with $k\geq2$ and minimum degree $\delta$,
then $\delta\geq2$.
\end{lemma}
\begin{proof}{ Suppose on the contrary that, there exists a vertex
$u\in V(G)$ with $deg(u)=1$. Let $v$ be the unique neighbor of $u$
and $T\subseteq V(G)$ be a subset of $V(G)$ with $|T|=k$ and
$u,v\in T$. Since $G$ is a randomly $k$-dimensional graph,
$T\setminus \{v\}$ is not a resolving set for $G$. Thus, there
exists a pair of vertices $x,y\in V(G)$ such that, $d(x,v)\neq
d(y,v)$ and $d(x,t)=d(y,t)$, for each $t\in T\setminus \{v\}$.
Hence, $d(x,u)=d(y,u)$. Clearly, if $u\in\{x,y\}$, then
$d(x,u)\neq d(y,u)$, which is a contradiction. Consequently,
$u\notin\{x,y\}$. Therefore, $d(x,u)=d(x,v)+1$ and
$d(y,u)=d(y,v)+1$. Thus, $d(x,v)=d(y,v)$. This  contradiction
implies that $\delta\geq2$. }\end{proof}
\begin{thm}\label{2connected} If $k\geq2$, then every randomly $k$-dimensional graph
is $2$-connected.
\end{thm}
\begin{proof}{Suppose on the contrary that $u$ is a cut vertex in
$G$. Let $G_1$ be a connected component of $G\setminus \{u\}$. Set
$H_2:=G\setminus V(G_1)$ and $H_1:=\langle
V(G_1)\cup\{u\}\rangle$, the induced subgraph by $V(G_1)\cup\{u\}$
of $G$. Note that, for each $x\in V(H_1)$ and each $y\in V(H_2)$,
$d(x,y)=d(x,u)+d(u,y)$. By Lemma~\ref{delta>2}, $G$ does not have
any vertex of degree $1$. Therefore, $|V(H_1)|\geq3$ and
$|V(H_2)|\geq3$. Let $a,b\in V(H_2)$ and $V(H_1)$ resolves
$\{a,b\}$. Then, there exists a vertex $w\in V(H_1)$ such that,
$d(a,w)\neq d(b,w)$. Thus, $d(a,u)+d(u,w)\neq d(b,u)+d(u,w)$, that
is $d(a,u)\neq d(b,u)$. Hence, $V(H_1)$ resolves a pair of
vertices of $V(H_2)$ if and only if $u$ resolves this pair. If
$V(H_1)$ is a resolving set for $G$, then $\{u\}$ is a resolving
set for $H_2$. Therefore, by Theorem~\ref{B=1,B=n-1}, $H_2$ is a
path. Since $|V(H_2)|\geq3$, $G$ has a vertex of degree $1$, which
contradicts Lemma~\ref{delta>2}. Hence, $\beta(H_2)\geq2$ and
$V(H_1)$ does not resolve $G$. Now,
one of the following two cases can be happened.\vspace{3mm}\\
1. $u$ belongs to a basis of $H_2$. In this case $u$ along with
$\beta(H_2)-1$ vertices of $V(H_2)\setminus \{u\}$ forms a basis
$T$ of $H_2$. Since $\beta(H_2)\geq2$, there exists a vertex $x\in
T\setminus\{u\}$. Note that, $T\cup V(H_1)\setminus\{x\}$ is not a
resolving set for $G$, otherwise $T\setminus\{x\}$ is a resolving
set for $H_2$ of size $\beta(H_2)-1$. Thus,
$$res(G)\geq |T\cup V(H_1)|=\beta(H_2)+|V(H_1)|-1.$$ Now, let
$z\in V(G_1)$. Since $|V(H_1)|\geq3$ and $G_1$ is a connected
component of $G\setminus\{u\}$, $z$ has a neighbor in $G_1$, say
$v$. Therefore, $d(z,v)=1\neq d(y,v)$ for each $y\in
V(H_2)\setminus\{u\}$. Hence, the set $T\cup V(H_1)\setminus\{z\}$
is a resolving set for $G$. Thus, $$\beta(G)\leq |T\cup
V(H_1)\setminus\{z\}|= \beta(H_2)+|V(H_1)|-2.$$ Consequently,
$\beta(G)<res(G)$, which is a contradiction.
\vspace{3mm}\\
2. $u$ does not belong to any basis of $H_2$. Let $T$ be a basis
of $G$ and $x\in T$. Therefore, $T\cup V(H_1)\setminus\{x\}$ is
not a resolving set for $G$. Hence,
$$res(G)\geq |T\cup V(H_1)|=\beta(H_2)+|V(H_1)|.$$ Now, let
$z\in V(G_1)$. Similar to the previous case, $T\cup
V(H_1)\setminus\{z\}$ is a resolving set for $G$. Thus,
$$\beta(G)\leq |T\cup V(H_1)\setminus\{z\}|=
\beta(H_2)+|V(H_1)|-1.$$ Therefore, $\beta(G)<res(G)$, which is
impossible.
\par Consequently, $G$ does not have any cut vertex.
}\end{proof}
\begin{thm}\label{deg2 not adjacent} If $G$ is a randomly $k$-dimensional
graph with $k\geq4$, then there are no adjacent vertices of degree
$2$ in $G$.
\end{thm}
\begin{proof}{ Suppose on the contrary that $G$ has adjacent
vertices of degree $2$. Therefore, there is an induced subgraph
$P_r=(a_1,a_2,\ldots,a_r),~r\geq2$, such that, for each $i$,
$1\leq i\leq r$, $deg(a_i)=2$ in $G$. Let $x,y\in V(G)\setminus
V(P_r)$ and $x\sim a_1$, $y\sim a_r$. Since $k\geq4$, $G$ is not a
cycle. Thus, Theorem~\ref{2connected} implies that $x\neq y$,
otherwise, $x=y$ is a cut vertex in $G$. By assumption, $G$ has a
basis $B=\{x,y,a_i,a_j\}\cup T$, where $1\leq i\neq j\leq r$ and
$T$ is a subset of $V(G)\setminus \{x,y,a_i,a_j\}$ with $|T|=k-4$.
Now, one of the following cases can be happened.
\vspace{3mm}\\
1. $r$ is odd. Let
$B_1=B\cup\{a_{r+1\over2}\}\setminus\{a_i,a_j\}$. We claim that,
$B_1$ is a resolving set for $G$. Otherwise, there exist vertices
$u,v\in V(G)$ with $r(u|B_1)=r(v|B_1)$. If $v\in V(P_r)$ and
$u\notin V(P_r)$, then $d(v,a_{r+1\over2})\leq {r-1\over2}$ and
$d(u,a_{r+1\over2})\geq {r+1\over2}$. Hence, $r(u|B_1)\neq
r(v|B_1)$, which is a contradiction. Therefore, both of vertices
$u$ and $v$ belong to $V(P_r)$ or $V(G)\setminus V(P_r)$. If
$u,v\in V(P_r)$, then, $d(u,a_{r+1\over2})=d(v,a_{r+1\over2})$
implies $u,v\in\{a_{{r+1\over2}-i},a_{{r+1\over2}+i}\}$ for some
$i$, $1\leq i\leq{r-1\over2}$. On the other hand,
$d(x,a_{{r+1\over2}-i})={r+1\over2}-i$ and
$d(x,a_{{r+1\over2}+i})=\min\{{r+1\over2}+i,{r+1\over2}-i+d(x,y)\}$.
If ${r+1\over2}+i\leq{r+1\over2}-i+d(x,y)$, then
$d(x,a_{{r+1\over2}-i})\neq d(x,a_{{r+1\over2}+i})$, which is a
contradiction. Thus, ${r+1\over2}-i+d(x,y)<{r+1\over2}+i$ and
hence, ${r+1\over2}-i+d(x,y)={r+1\over2}-i$, because
$d(x,a_{{r+1\over2}-i})=d(x,a_{{r+1\over2}+i})$. Therefore,
$d(x,y)=0$, which  contradicts $x\neq y$. Thus, $u,v\in
V(G)\setminus V(P_r)$. Since $r(u|B_1)=r(v|B_1)$ and $B$ is a
resolving set for $G$, there exists a vertex in $B\setminus
B_1=\{a_i,a_j\}\setminus\{a_{r+1\over2}\}$ which resolves
$\{u,v\}$. By symmetry, we can assume $a_i$ resolves $\{u,v\}$.
Therefore, $d(u,a_i)\neq d(v,a_i)$, $d(u,x)=d(v,x)$, and
$d(u,y)=d(v,y)$. But,
$$d(u,a_i)=\min\{d(u,x)+d(x,a_i),d(u,y)+d(y,a_i)\},$$and
$$d(v,a_i)=\min\{d(v,x)+d(x,a_i),d(v,y)+d(y,a_i)\}.$$If
$d(u,x)+d(x,a_i)\leq d(u,y)+d(y,a_i)$ and $d(v,x)+d(x,a_i)\leq
d(v,y)+d(y,a_i)$, then $d(u,x)+d(x,a_i)\neq d(v,x)+d(x,a_i)$,
which implies $d(u,x)\neq d(v,x)$, a contradiction. Similarly, if
$d(u,y)+d(y,a_i)\leq d(u,x)+d(x,a_i)$ and $d(v,y)+d(y,a_i)\leq
d(v,x)+d(x,a_i)$, then $d(u,y)\neq d(v,y)$, which is a
contradiction. Therefore, by symmetry, we can assume
$d(u,x)+d(x,a_i)\leq d(u,y)+d(y,a_i)$ and $d(v,y)+d(y,a_i)\leq
d(v,x)+d(x,a_i)$.
Thus,$$d(u,a_i)=d(u,x)+d(x,a_i)=d(v,x)+d(x,a_i)\geq d(v,a_i),$$
and$$d(v,a_i)=d(v,y)+d(y,a_i)=d(u,y)+d(y,a_i)\geq d(u,a_i).$$
These imply that $d(u,a_i)=d(v,a_i)$, which is a contradiction.
Therefore, $B_1$ is a resolving set for $G$ with cardinality
$k-1$.\vspace{3mm}\\
2. $r$ is even. Let
$B_2=B\cup\{a_{r\over2}\}\setminus\{a_i,a_j\}$. Similar to the
previous case, $B_2$ is a resolving set for $G$ with cardinality
$k-1$.
\par In both cases, we get a contradiction to the assumption
that $G$ is a randomly $k$-dimensional graph. Therefore, there are
no adjacent vertices of degree $2$ in $G$. }\end{proof}
\begin{thm}\label{|T|=k-1} If $G$ is a randomly $k$-dimensional
graph and $T$ is a separating set of $G$ with $|T|=k-1$, then
$G\setminus T$ has exactly two connected components and for each
pair of vertices $u,v\in V(G)\setminus T$ with $r(u|T)=r(v|T)$,
$u$ and $v$ belong to different components.
\end{thm}
\begin{proof}{ Since $\beta(G)=k$ and $|T|=k-1$, there exist two
vertices $u,v\in V(G)\setminus T$ with $r(u|T)=r(v|T)$. Let $H$ be
a connected component of $G\setminus T$ for which $u\notin H$ and
$v\notin H$. If $w\in H$, then there exist two vertices $s,t\in T$
such that, $d(u,w)=d(u,s)+d(s,w)$ and $d(v,w)=d(v,t)+d(t,w)$.
Since $r(u|T)=r(v|T)$, we have $d(u,s)=d(v,s)$ and
$d(u,t)=d(v,t)$. Therefore,
$$d(u,w)=d(u,s)+d(s,w)=d(v,s)+d(s,w)\geq d(v,w).$$ And
$$d(v,w)=d(v,t)+d(t,w)=d(u,t)+d(t,w)\geq d(u,w).$$ Hence,
$d(u,w)=d(v,w)$. Thus, $r(u|T\cup\{w\})=r(v|T\cup\{w\})$.
Consequently, $T\cup\{w\}$ is not a resolving set for $G$ and
$|T\cup\{w\}|=k$. This contradicts the assumption that $G$ is
randomly $k$-dimensional. Therefore, $G\setminus T$ has exactly
two components and $u$ and $v$ belong to different components.
}\end{proof}
\begin{cor}\label{Delta>k}
If $G$ is a randomly $k$-dimensional graph with $k\geq2$, then
$\Delta(G)\geq k$.
\end{cor}
\begin{proof}{ If $G=K_n$, then $\Delta(G)=n-1=k$. Now let $G\neq
K_n$. Suppose on the contrary that $\Delta(G)\leq k-1$. Let $u\in
V(G)$, $deg(u)=\Delta(G)$, and $T$ be a subset of $V(G)$ with
$|T|=k-1$ and $N(u)\subseteq T$. By Theorem~\ref{|T|=k-1},
$G\setminus T$ has exactly two connected components, of which one
of them is $\{u\}$. Since $|T|=k-1$ and $\beta(G)=k$, there exist
two vertices $x,y\in V(G)\setminus T$ such that, $r(x|T)=r(y|T)$.
By Theorem~\ref{|T|=k-1}, $x$ and $y$ belong to different
components. Therefore, one of them is $u$, say $x=u$. Since
$r(u|T)=r(y|T)$, we have $N(u)\subseteq N(y)$. By
Corollary~\ref{not twin}, $G$ does not have any pair of vertices
$u,v$  with $N(u)\setminus\{v\}=N(v)\setminus\{u\}$. Hence,
$N(u)\subset N(y)$, this contradicts $deg(u)=\Delta(G)$.
Therefore, $\Delta(G)\geq k$. }\end{proof}
\begin{cor}\label{deg not adjacent>k}
If $u$ and $v$ are two non-adjacent vertices in  a randomly
$k$-dimensional graph, then  $deg(u)+deg(v)\geq k$.
\end{cor}
\begin{proof}{ If $|N(u)\cup N(v)|\leq k-1$, then let $T$ be a subset of
$V(G)\setminus\{u,v\}$ with $|T|=k-1$ and $N(u)\cup N(v)\subseteq
T$. By Theorem~\ref{|T|=k-1}, $G\setminus T$ has exactly two
connected components $\{u\}$ and $\{v\}$. Hence, $|T|=n-2$. This
implies that $k=n-1$ and by Theorem~\ref{res=1,res=n-1}, $G=K_n$.
Consequently, $u\sim v$, which is a contradiction. Thus,
$deg(u)+deg(v)\geq|N(u)\cup N(v)|\geq k$.
 }\end{proof}
 \begin{thm}\label{omega<k+1} If $G$ is a randomly $k$-dimensional
 graph of order at least $2$, then $\omega(G)\leq k+1$. Moreover, $\omega(G)=k+1$ if and
 only if $G=K_n$.
 \end{thm}
 \begin{proof}{ Let $H$ be a clique of size $\omega(G)$ in $G$ and
 $T$ be a subset of $V(H)$ with $|T|=\omega(G)-2$. If
 $T=V(H)\setminus\{u,v\}$, then $r(u|T)=(1,1,\ldots,1)=r(v|T)$.
 Therefore, $T$ is not a resolving set for $G$. Since $G$ is a
 randomly $k$-dimensional graph, $|T|\leq k-1$. Thus,
 $\omega(G)-2=|T|\leq k-1$. Consequently, $\omega(G)\leq k+1$.
 \par Clearly, if $G=K_n$, then $\omega(G)=k+1$. Conversely, let
 $\omega(G)=k+1$. If $G\neq K_n$, then there exists  a vertex $x\in
 V(G)\setminus V(H)$ such that, $x$ is adjacent to some vertices
 of $V(H)$, because $G$ is connected. Since $|V(H)|=\omega(G)$,
 $x$ is not adjacent to all vertices of $V(H)$. If there exist
 vertices $y,z\in V(H)$ such that, $y\nsim x$ and $z\nsim x$, then
 $d(x,y)=d(x,z)=2$, because
 $x$ is adjacent to some vertices of $H$. Let
 $S=\{x\}\cup V(H)\setminus\{y,z\}$. Therefore,
 $r(y|S)=(2,1,1,\ldots,1)=r(z|S)$. Thus, $S$ is not a resolving
 set for $G$ and $|S|=k$, which is  a contradiction. Hence, $x$ is
 adjacent to $\omega(G)-1$ vertices of $H$.
\par On the other
 hand, $x$ is adjacent to at most one vertex of $H$. Otherwise, there exist
 vertices $s,t\in V(H)$ such that, $s\sim x$ and $t\sim x$. Let
 $R=\{x\}\cup V(H)\setminus\{s,t\}$. Therefore,
 $r(s|R)=(1,1,\ldots,1)=r(t|R)$. Thus, $R$ is not a resolving
 set for $G$ and $|R|=k$, which is  a contradiction. Consequently,
 $\omega(G)=2$ and  $k=\omega(G)-1=1$. Therefore, $G=K_2$,
 which contradicts $G\neq K_n$. Hence, $G=K_n$. }\end{proof}
\begin{lemma}\label{at most k-1 common neighbor} If $res(G)=k$,
 then each two vertices of $G$
 have at most $k-1$ common neighbors.
\end{lemma}
\begin{proof}{ Let $u,v\in V(G)$ and $T=N(u)\cap N(v)$. Thus,
$r(u|T)=(1,1,\ldots,1)=r(v|T)$. Therefore, $T$ is not a resolving
set for $G$. Since $G$ is a randomly $k$-dimensional graph,
$|N(u)\cap N(v)|=|T|\leq k-1$.
 }\end{proof}
 \begin{thm}\label{delta<n-2} If $G\neq K_n$ is a randomly $k$-dimensional
 graph of order $n$, then $\Delta(G)\leq n-2$.
 \end{thm}
 \begin{proof}{ Suppose on the contrary that there exists a vertex
 $u\in V(G)$ with $deg(u)=n-1$. For each $T\subseteq
 V(G)\setminus \{u\}$ with $|T|=k-1$, the set $T\cup\{u\}$ is a
 resolving set for $G$ while, $T$ is not a
 resolving set for $G$. Hence, there exist vertices $x,y\in
 V(G)\setminus T$ such that, $r(x|T)=r(y|T)$ and $d(x,u)\neq
 d(y,u)$. Since $u$ is adjacent to all vertices of $G$, we have
 $u\in\{x,y\}$, say $x=u$. Thus, $r(y|T)=r(u|T)=(1,1,\ldots,1)$.
 By Lemma~\ref{at most k-1 common neighbor}, $|N(u)\cap N(y)|\leq
 k-1$. Hence, $deg(y)\leq k$, because $u$ is adjacent to all vertices of $G$.
This gives, $N(y)= T\cup\{u\}$.
 \par Now, let
 $S=T\cup\{y\}\setminus\{v\}$, for an arbitrary vertex $v\in T$.
 Since $|S|=k-1$, $S$ is not a resolving set for $G$. Therefore,
 there exist vertices $a,b\in V(G)\setminus S$ such that,
 $r(a|S)=r(b|S)$. Since $S\cup\{u\}$ is a resolving set for $G$,
 we have $d(a,u)\neq d(b,u)$. Hence, $u\in\{a,b\}$, say $b=u$.
 Thus, $r(a|S)=r(u|S)=(1,1,\ldots,1)$. Consequently, $a\sim y$.
  Therefore, $a\in T$, because $N(y)= T\cup\{u\}$
 and $a\neq u$. Hence, $a\in(V(G)\setminus S)\cap T=\{v\}$, that is $a=v$.
 Thus, $v$ is adjacent to all vertices of $T\setminus\{v\}$. Since $v$ is an
 arbitrary vertex of $T$, $T$ is a clique. Therefore,
$T\cup\{u,y\}$ is a clique of size $k+1$ in $G$. Consequently, by
Theorem~\ref{omega<k+1}, $G=K_n$, which is a contradiction. Thus,
 $\Delta(G)\leq n-2$.
}\end{proof}


\begin{thebibliography}{99}
\bibitem{net2}
{\it Z. Beerliova, F. Eberhard, T. Erlebach, A. Hall, M. Hoffmann,
M. Mihal'ak and L.S. Ram}, Network dicovery and verification,
 {\it IEEE Journal On Selected Areas in Communications} {\bf24(12)} (2006) 2168-2181.
\bibitem{cartesian product} {\it  J. Caceres, C. Hernando, M. Mora, I.M. Pelayo,
M.L. Puertas, C. Seara and D.R. Wood}, {\it On the metric
dimension of cartesian products of graphs},~SIAM Journal on
Discrete Mathematics {\bf 21(2)} (2007) 423-441.
\bibitem{bounds}
{\it G.G. Chappell, J. Gimbel and C. Hartman}, {\it Bounds on the
metric and partition dimensions of a graph}, Ars Combinatorics
{\bf 88} (2008) 349-366.
\bibitem{Ollerman} {\it  G. Chartrand, L. Eroh, M.A. Johnson and O.R. Ollerman},
 {\it Resolvability in graphs and the metric dimension of a graph},~Discrete Applied Mathematics {\bf 105} (2000) 99-113.
\bibitem{sur1}
{\it G. Chartrand and P. Zhang}, {\it The theory and applications
of resolvability in graphs.} A survey. In Proc. 34th Southeastern
International Conf. on Combinatorics, Graph Theory and Computing
 {\bf 160} (2003) 47-68.
\bibitem{basis}
{\it G. Chartrand and P. Zhang}, {\it On the chromatic dimension
of a graph}, Congressus Numerantium {\bf 145} (2000) 97-108.
\bibitem{res(G)}
{\it G. Chartrand, C. Poisson and P. Zhang}, {\it Resolvability
and the upper dimension of graphs}, Computers and Mathematics with
Applications {\bf 39} (2000) 19-28.
\bibitem{Harary}{\it F. Harary and R.A Melter}, {\it On the metric dimension of a graph},~Ars Combinatorics {\bf 2} (1976) 191-195.
\bibitem{extermal} {\it C. Hernando, M. Mora, I.M. Pelayo, C. Seara and D.R. Wood},
{\it Extemal graph theory for metric dimension and diameter},~ The
Electronic Journal of Combinatorics {\bf } (2010) \#R30.

\bibitem{landmarks}{\it S. Khuller, B. Raghavachari and A. Rosenfeld},
 {\it Landmarks in graphs},~Discrete Applied Mathematics
{\bf 70(3)} (1996) 217-229.
\bibitem{digital}{\it R.A. Melter and I. Tomescu},
 {\it Metric bases in digital geometry}, Computer Vision Graphics and
 Image Processing
{\bf 25} (1984) 113-121.
\bibitem{sur2} {\it V. Saenpholphat and P. Zhang}, {\it Conditional resolvability in graphs:
a survey}, International Journal of Mathematics and Mathematical
Sciences {\bf 38} (2004) 1997-2017.
\bibitem{coin}{\it A. Sebo and E. Tannier},
 {\it On metric generators of graphs}, Mathematics of Operations
 Research
{\bf 29(2)} (2004) 383-393.
\bibitem{Slater1975} {\it P.J. Slater}, {\it Leaves of trees},
Congressus Numerantium {\bf 14} (1975) 549-559.
\bibitem{west} {\it D.B. West}, Introduction to graph theory,
{\it Prentice Hall Inc. Upper Saddle River, NJ 07458, Second
Edition} (2001).
\end{thebibliography}
\end{document}